\documentclass[12pt]{article}
\usepackage{amssymb}
\usepackage{amsmath}
\usepackage{amstext}
\usepackage{amsfonts}
\usepackage{amsthm}
\usepackage{amscd}
\numberwithin{equation}{section}

\swapnumbers
\theoremstyle{plain}
\newtheorem{thm}{Theorem}[section]
\newtheorem{lem}[thm]{Lemma}
\newtheorem{prop}[thm]{Proposition}
\newtheorem{cor}[thm]{Corollary}
\newtheorem{definition}[thm]{Definition}
\newtheorem{conj}[thm]{Conjecture}
\theoremstyle{remark}
\newtheorem{rem}[thm]{Remark}
\newtheorem{que}[thm]{Question}

\def\mod0{{\mathrm M}_X}

\def\cO{\mathcal{O}}
\def\g{\mathfrak{g}}

\def\map#1{\ \smash{\mathop{\longrightarrow}\limits^{#1}}\ }

\def\LL{\mathcal{L}}
\def\XX{\mathcal{X}}
\def\MM{\mathrm{M}}
\def\NN{\mathrm{N}}
\def\GG{\mathbb{G}}
\def\PP{\mathbb{P}}

\def\ZZ{\mathbb{Z}}

\def\QQ{\mathbb{Q}}

\def\DD{\mathcal{D}}
\def\L{\mathbb{L}}
\def\Gr{\mathrm{Gr}}
\def\De{\Delta}
\def\T{\mathbf{T}}

\def\div{\mathrm{div}}
\def\deg{\mathrm{deg} \:}

\def\sym{\mathrm{Sym}}
\def\ker{\mathrm{ker} \:}

\def\kum{\mathrm{Kum} \:}
\def\dim{\mathrm{dim} \:}

\def\im{\mathrm{im} \:}

\def\Hom{\mathrm{Hom}}
\def\supp{\mathrm{supp} \:}
\def\Spec{\mathrm{Spec}}

\def\lra{\longrightarrow}

\def\ra{\rightarrow}

\def\lms{\longmapsto}
\def\isom{\map{\sim}}
\def\kum{\mathrm{Kum}}
\def\GL{\mathrm{GL}}

\def\pp{\mathbb{P}}

\def\M{\mathcal{M}}

\title{The action of the Frobenius map on rank $2$ vector bundles
in characteristic $2$}
\author{Y. Laszlo and C. Pauly}
\begin{document}
\maketitle

\section{Introduction}
Let $X$ be a smooth algebraic curve of genus $g$ over a field $k$
of characteristic $p>0$. The behaviour of semi-stable bundles
with respect to the absolute Frobenius $F_a$ remains mysterious
if $g \geq 2$. Let us briefly explain why this question should be
of interest. Start with a continuous representation $\rho$ of the
algebraic fundamental group in $\GL_r(\bar{k})$, where $\bar{k}$
is the algebraic closure of $k$. Let $E_\rho$ be the
corresponding rank $r$ bundle over $X$. Then all the bundles
$F_a^{(n)*}E_\rho$, for $n > 0$, where $F_a^{(n)}$ denotes the
$n$-fold composite $F_a \circ \cdots \circ F_a$, are semi-stable.
Conversely, assuming that $k$ is finite, let $E$ be a semi-stable
rank $r$ bundle defined over $\bar{k}$. Because the set of
isomorphism classes of semi-stable bundles of degree $0$ over
$X_{k'}$, where $k'$ is any finite extension of $k$, is finite,
one observes (see \cite{ls}) that some twist of $E$ comes from a
representation as above.  Therefore, if one is interested in
unramified continuous representations of the Galois group over $k$
of a global field $k(X)$ in characteristic $p$, it is natural to
look at Frobenius semi-stable bundles, that is those whose
pull-backs by $F_a^{(n)}$ are all semi-stable. This condition is
stable by tensor product (\cite{mi} section 5), which is not
usually the case for ordinary semi-stability in positive
characteristic.

\bigskip

Assume that $k$ is arbitrary of characteristic $p$. Let us
emphasize that the locus of Frobenius semi-stable bundles
of degree $0$ is a countable intersection of open subsets
of the coarse moduli scheme of semi-stable vector bundles of
degree $0$ over $X$. In particular, it is not clear at all
if the Frobenius semi-stable points are dense in general.
This could a priori depend on the arithmetic of the base
field $k$. We would like to study the dynamics of the
Frobenius, namely the sequence $n \lms F_a^{(n)*}E$ for
a given vector bundle $E$ over $X$.

\bigskip

If $k$ is a discrete valuation field of characteristic $p$ and if
$X$ is a Mumford curve, the situation is well understood
\cite{fal}, \cite{gie}: among other results, it is shown for
arbitrary genus and characteristic that there exist semi-stable
bundles which are destabilized by the Frobenius $F_a$, that $E
\lms F^*_a E$ induces a surjective (rational) map on the moduli
space of semi-stable rank $r$ vector bundles of degree $0$, and
that the set of bundles coming from continuous representations of
the algebraic fundamental group is dense. Another case, which is
studied in the literature, are elliptic curves (see e.g.
\cite{oda},\cite{nsb}). For example, it can easily be shown that
over an elliptic curve semi-stability is preserved under
pull-back by Frobenius and that a stable bundle of rank $r$ and
degree $d$ is Frobenius stable if and only if $pd$ and $r$ are
coprime. But in general, not much seems to be known.

\bigskip

In this paper we study the action of the Frobenius map $F_a$ in a
very particular case: $X$ is an ordinary curve of genus $2$
defined over an  algebraically closed field $k$ of characteristic
$2$. In that case the coarse moduli space $\mod0$ of semi-stable
rank $2$ vector bundles of trivial determinant is known to be the
projective space $\pp^3$ and we show (Proposition
\ref{computations}) that the Frobenius map identifies with a
rational map $\pp^3 \dashrightarrow \pp^3$ given by quadratic
polynomials, which are explicitly computed in terms of generalized
theta constants of the curve $X$. This result allows us to deduce
that, if the determinant is fixed, the Frobenius map is not
surjective and that the set of Frobenius semi-stable bundles is
Zariski dense in $\pp^3$ (Proposition \ref{dynamics}). On the
other hand, the Frobenius map becomes surjective if the
determinant is not fixed (Proposition \ref{surjectivity}).

\bigskip

As a by-product of our computations, we also get the explicit
equation of Kummer's quartic surface (section 4) and
a description of the Frobenius map acting on the moduli space
of rank $2$ vector bundles of fixed odd degree determinant
(section 7).

\bigskip

For other aspects of this problem see the recent
article \cite{xia}.

\bigskip

We would like to thank M.S. Narasimhan and C.S. Seshadri
who explained to the authors why the classical identification
$\mod0  = \pp^3$ of \cite{NarRam} remains valid in characteristic $2$. We
especially thank M. Raynaud for having introduced us to these
questions and for his encouragements.

\section{Review of Theta groups and
Theta divisors in characteristic $p$}
Let $k$ be a perfect field of characteristic $p>0$.
By schemes we implicitly mean $k$-scheme.
\subsection{Relative Frobenius}

Let $A$ be an ordinary abelian variety of dimension $g$. We shall
denote by $F$ the relative Frobenius morphism
$$F : A \lra A_1$$
which is a purely inseparable $k$-isogeny of degree $p^g$.
Its kernel $\hat{G}$ is a subgroup scheme of the group scheme
$A[p]$ of $p$-torsion points of $A$ and $F$ identifies to the
quotient morphism $A \lra A/\hat{G}$. Because
$$\ker F =\hat{G}\subset A[p]=\ker [p],$$ the diagram
$$ A \lra A/\hat{G} \lra A/A[p]$$ becomes

$$ [p]: A \map{F} A_1 \map{V} A.$$

\begin{rem}
In the case where $A$ is a Jacobian, the map $V$, called Verschiebung,
is simply the pull-back by the relative Frobenius $F: X\ra X_1$.
\end{rem}

Let $G$ be the reduced part of $A[p]$. Because $A$ is ordinary, $G$
intersects $\hat{G}$ at the origin and one gets a canonical
decomposition
\begin{equation} \label{deckl}
A[p]= G \times \hat{G}
\end{equation}
and the relative Frobenius induces an isomorphism
$$F:\ \ker V =A[p]/\hat{G}\isom G.$$

\subsection{Theta group scheme}
We will use the basic results (and notations) of
Theta groups associated to line bundles on abelian varieties as in
\cite{mum2}. We assume that $A_1$ is
principally polarized and that $k$ is big
enough in order to have an isomorphism $G \cong (\ZZ/p\ZZ)^g$
(however, we do not choose such an isomorphism). Let $\Theta_1$ be
a (not necessarily symmetric) Theta divisor representing the
polarization. We denote by $M_1$ (resp. $L_1$) the line bundle
$\cO(\Theta_1)$ (resp. $\cO(p\Theta_1)$).
We consider the Theta group $G(L_1)$ associated to the line
bundle $L_1$ (see \cite{mum2}, page 221).
This group scheme is a central extension
$$ 0 \lra \GG_m \lra G(L_1) \lra K(L_1) \lra 0 $$
and $K(L_1)=A_1[p]$. We denote by $e^{L_1}: K(L_1) \times
K(L_1) \ra \GG_m$ the skew-symmetric form on $K(L_1)$ given by
the commutator taken in the Theta group $G(L_1)$.

\begin{lem}\label{L1}
There exists a unique ample line bundle $M$ on $A$ such that
$V^*M=L_1.$
\end{lem}

\begin{proof}
Because $K(L_1)$ is self-dual and contains $\hat{G}$, it is the
whole $A[p]$. One has to show (theorem 2, page 231 of \cite{mum2})
that the restriction of $G(L_1)$ to $G=\ker V$ is
 split (the uniqueness is obvious). Since $G$ is reduced, it is
 enough to define the splitting at the
level of $k$-points. Let $g \in G(k)$. Since $\mu_p(k) = \{1\}$ and
$k^*$ is divisible, there exists a unique $\tilde{g} \in G(L_1)$
over $g$ satisfying $\tilde{g}^p = 1$. Since $\mu_p(k) = \{1\}$
again, the restriction of $e^{L_1}$ to $G(k)\times G(k)$ is trivial
and therefore $G$ is isotropic for $e^{L_1}$. This implies that any
two elements $\tilde{g},\tilde{h}$ commute and the map $g
\mapsto \tilde{g}$ is a morphism of group schemes.
\end{proof}

The bundle $M$ defines a principal polarization on $A$ and the
$G$-invariant morphism
$$V^*M^p=V^*[(V_*L_1)^G]^{\otimes p}\lra L_1^p=M_1^{p^2}=V^*F^*M_1$$
defines an isomorphism
\begin{equation}\label{isoMFM}
M^p\isom F^*M_1.
\end{equation}

We denote by $G(L)$ the Theta group scheme associated to the line
bundle $L:= F^*M_1 = M^p$.

\begin{lem} \label{liftGdual}
The restriction of $G(L)$ to both $G$ and $\hat{G}$ is canonically
split.
\end{lem}

\begin{proof}
Because $L$ is the pull-back of $M_1$ by $F$, the restriction of
$G(L)$ to both $\hat{G}$ splits (see Theorem 2, page 231 of
\cite{mum2} again). Observe however that the splitting is
determined not only by $L$ but by $M_1$. For the splitting over
$G$, proceed as in lemma \ref{L1}.
\end{proof}

Therefore, the decomposition $A[p]=K(L)=\hat{G}\times G$ of
(\ref{deckl}) is symplectic in the sense that both $G$ and
$\hat{G}$ are isotropic for $e^L$. In particular, $e^L$ identifies
$\hat{G}$ and the Cartier dual of $G$.

\begin{rem}\label{Katz}
One can interpret more geometrically the action of $\hat{G}$ as
follows. By \cite{kat}, the line bundle $L$ comes with a canonical
$p$-integrable connection $\nabla$. The Lie algebra $\g$ of
$\hat{G}$ is the whole tangent space $\T_0A$. The link between the
infinitesimal action of $\g$ on $L$ and $\nabla$ is simply given by
$$v.l=\nabla_vl$$
where $l$ is a local section of $L$, $v$ a point of $\g$ and
$\nabla_vl$ is the $\nabla$-derivative of $l$ with respect to the
invariant vector field defined by $v$.
\end{rem}

Let $\theta$ be a non-zero global section of $M$. By the remark above,
since the section $\theta^p$ is annihilated by $\nabla$, it defines a
$\hat{G}$-invariant section of $F^*\Theta_1$, namely a non-zero
section $\theta_1$ of $M_1$ such that
\begin{equation}\label{theta=Ftheta}
 \theta^p=F^*\theta_1.
\end{equation}

The main theorem on the group scheme $G(L)$ and its representations
is the following structure theorem due to Mumford (for the
characteristic $p$ version we need, we refer to \cite{sek}):

\begin{thm}\label{Heisenberg}
The space of global sections $H^0(A,L)$ is the unique irreducible
representation of weight $1$ (i.e. $\GG_m$ acts by its natural
character) of the Theta group scheme $G(L)$.
\end{thm}

\subsection{Canonical Theta basis of $|p\Theta|$ and $|p\Theta_1|$}

In \cite{mum1} Mumford constructs canonical bases for any linear
system $\pp H^0(A,L)$ where $L$ is a line bundle of separable type.
Because of Theorem \ref{Heisenberg}, we can adapt his construction to line
bundles $L$ not of separable type. In our situation, namely $L =
\cO_A(p\Theta)$, where $\Theta$ is not necessarily symmetric,
we get the following
\begin{lem} \label{thetabasis}

(i) There exists a basis $\{X_g\}_{g \in G}$ of $H^0(A,L)$, unique up
to a multiplicative scalar, which satisfies the following relations
\begin{equation} \label{heiscom}
a.X_g = X_{g+a} \qquad \alpha.X_g = e^L(\alpha,g) X_g
\qquad \forall a,g \in G, \ \forall \alpha \in \hat{G}
\end{equation}

(ii) For every $g\in G$, there exists a unique section $Y_g\in
H^0(A_1,L_1)$ such that $X_g^p=F^*Y_g$. The family $\{Y_g\}_{g
\in G}$ is a basis of $H^0(A_1,L_1)$.
\end{lem}

\begin{proof}
Let us construct geometrically the basis $\{X_g\}_{g \in G}$. We
define $X_g$ by
$$X_g=g.F^*\theta_1=g.\theta^p.$$
By construction, one has $$a.X_g = X_{g+a} \qquad \forall a,g \in
G.$$ Because $X_0$ comes from $A_1=A/\hat{G}$, it is $\hat{G}$
invariant. The relations (\ref{heiscom}) follow. Because $H^0(A,L)$
is irreducible, $\{X_g\}_{g\in G}$ is a basis. If we have another
such basis $X'_g$, the endomorphism given by $X_g\longmapsto X'_g$
is $G(L)$-equivariant and therefore a scalar by Schur's lemma,
which proves (i).

The sections $X_g^p$ live in $H^0(A,L^p)=H^0(A,F^*L_1)^{\hat{G}}$
by (i) and (ii) follows.
\end{proof}

We will identify $H^0(A,L)$ and $H^0(A_1,L_1)$ with their duals
(more precisely $\pp H^0(A,L)$ and $\pp H^0(A,L)$) using these bases.

\begin{rem} \label{thetabasispowers}
More generally, we can construct Theta bases for any power
$L = \cO(p^l\Theta)$ with $l \geq 1$. We denote by
$G_l$ the reduced part of $\ker [p^l]$. Note that we have
$G_l = (\ZZ/p^l\ZZ)^g$, $G_1 = G$ and, for any $l \geq 1$,
we have an exact sequence $0 \lra G_1 \lra G_{l+1} \map{p}
G_l \lra 0$. As above, we prove that there exist canonical
bases $\{ X_g^{(l)} \}_{g \in G_l}$ of $H^0(A,L)$ and
 $\{ Y_g^{(l)} \}_{g \in G_l}$ of $H^0(A_1,L_1)$ satisfying
relations \eqref{heiscom}.
\end{rem}

\subsection{Addition formula in characteristic $2$}

From now we assume that $M$ is symmetric and that $p=2$. We
will explain how to obtain an addition formula in this
context. The method essentially goes as in the proof of Lemma
1.2 of \cite{sek}, which is a generalization of Mumford's
arguments.

Consider the homomorphism
\begin{equation} \label{morxi}
\xi:\
\left\{
\begin{array}{ccc}
    A\times A & \lra &A\times A \\
    (x,y)&\lms& (x-y,x+y)
  \end{array}
\right.
\end{equation}

(in the notation of {\it loc. cit.}, $a=b=1$). By the See-Saw
theorem, we have an isomorphism
\begin{equation} \label{isoxi}
\xi^*(M\boxtimes M) \cong M^2\boxtimes M^2
\end{equation}
hence we get an injection
\begin{equation*} \label{injxi}
\xi^*: H^0(A,M) \otimes H^0(A,M) \hookrightarrow
H^0(A,M^2) \otimes H^0(A,M^{2})
\end{equation*}

Let $\theta$ be the unique (up to a scalar) section of $H^0(A,M)$
and consider the Theta basis $\{ X_g \}_{g \in G}$ \eqref{heiscom}
of $H^0(A,L) = H^0(A,M^2)$. We have an exact sequence
$$ 1 \lra \GG_m \lra G(M^2 \boxtimes M^{2}) \lra
A[2] \times A[2] \lra 1 .$$
The kernel $K = \ker \xi = A[2]$ sits
diagonally in $A[2] \times A[2]$ and the isomorphism \eqref{isoxi}
determines a lift $K^*
\subset G(M^2
\boxtimes M^{2})$ of $K$.

The main point is that $K=A[2]$ has a G\"opel system in classical
terminology of theta functions, namely
$$K \cong G \times \hat{G}$$
where both $G$ and $\hat{G}$ are isotropic.
In {\it loc. cit.}, this existence follows from the condition
$p\nmid a+b$, which is the only reason to put this arithmetic
condition, which of course is not fulfilled here. Via the
projection map
$$ G(M^2) \times G(M^{2}) \lra G(M^2 \boxtimes M^{2}) $$
the lift $K^*$ induces lifts of $G$ and $\hat{G}$ into $G(M^{2})$.
Because the $\hat{G}$-invariant part of $H^0(A,M^2)$ is generated
by $X_0$, Sekiguchi's result gives

\begin{lem}[Sekiguchi] \label{addform2}
Normalizing $\xi$ suitably, one has the formula
\begin{equation}\label{addition}
\xi^*(\theta^{} \boxtimes \theta) =
 \sum_{h \in G} X_h \boxtimes X_h
\end{equation}
\end{lem}
In other words, we have
$$\theta(x-y)\theta(x+y)=\sum_{g\in G}X_g(x)X_g(y),\ x,y\in A.$$

Let us define the Kummer morphism $\kum_A:$ $A \ra |2\Theta|$,
$y \mapsto T^*_y \Theta + T^*_{-y} \Theta$, where $T_y$ denotes
translation by $y$. Then we can write

\begin{equation*}
  \div (\sum_gX_g(y)X_g )=T^*_y\Theta+T^*_{-y}\Theta=\kum_A(y), y\in A.
\end{equation*}
Using the relation $F^*\Theta_1=2\Theta$, one gets the analogous
relation
\begin{equation*}
  \div (\sum_gY_g(y)Y_g)=T^*_y\Theta_1+T^*_{-y}\Theta_1=
\kum_{A_1}(y), y\in A_1.
\end{equation*}

The element \eqref{addition} induces a linear isomorphism
$H^0(A,M^2)^* \map{\sim} H^0(A,M^2)$, which allows us to identify
both spaces.

\begin{cor}\label{kummer}
With the identification above, the complete
linear system $\varphi_L$ (resp. $\varphi_{L_1}$) is the Kummer morphism
$\kum_A$ (resp. $\kum_{A_1}$).
\end{cor}

\subsection{The Theta divisor $\Theta$ in characteristic $2$}
\label{thetadiv}
From now we assume that $p=2$ and that $X$ is an ordinary
curve of genus $g$, whose Jacobian is denoted by $J$. Let $B$
be the theta-characteristic of $X_1$ defined by the exact sequence
(\cite{Ray} section 4)
\begin{equation} \label{defB}
0 \lra \cO_{X_1} \lra F_* \cO_X \lra B \lra 0
\end{equation}
and we denote by $\Theta_1 \subset J_1$ the symmetric Theta
divisor determined by $B$, i.e.,
\begin{equation} \label{defthetaB}
\supp \Theta_1 = \{ N \in J_1 : \ h^0(X_1, B \otimes N) > 0 \}.
\end{equation}
We denote by $\Theta$ the Theta divisor on $J$ obtained from
$\Theta_1$ (Lemma \ref{L1}). Note that we have $\Theta_1 = \iota^*
\Theta$, where $\iota: J_1 \ra J$ is the $k$-semi-linear isomorphism.
Let $R$ be the ring of dual numbers $k[\epsilon]$ with $\epsilon^2 = 0$.
We recall that the Lie algebra $\g = \hat{G}(R)$ identifies with
the tangent space $\T_0J$. For any tangent vector $v \in \g$ and
$g \in G$ we still denote by $e^L(v,g) \in k$
the derivative of
$$ \hat{G} \lra \GG_m, \qquad \bar{g} \lms e^L(\bar{g},g). $$
Writing $\bar{g} = 1 + \epsilon v$, with $v \in \g$, we obtain
$\forall v \in \g$,$\forall g,h \in G$,
$e^L(v,g+h) = e^L(v,g) + e^L(v,h)$. We recall some well-known
facts about $\Theta$.
\begin{lem} \label{thetainf}
\begin{itemize}
\item[(i)]
The Theta divisor $\Theta_1$ (and $\Theta$) passes through any
non-zero $g \in G$.
\item[(ii)]
A non-zero $g$ is a smooth point of $\Theta_1$ (and $\Theta$) if
and only if $h^0(X_1,B \otimes g) =1$.
\item[(iii)]
Assuming (ii), the tangent space at $g$ $\T_g \Theta_1 \subset \T_g J_1$
is defined by the linear form $e^L(.,g)$ on $\g$. Alternatively,
identifying the tangent space $\T_g J$ with $H^1(X,\cO_X)$, the
tangent space $\T_g \Theta$ is the image of the injective $k$-linear
map induced by the relative Frobenius,i.e.,
$$ \T_g \Theta = \im \ F:  H^1(X_1,g) \hookrightarrow H^1(X,\cO_X).$$
\end{itemize}
\end{lem}

\begin{proof}
Part (i) follows immediately from \eqref{defB} and \eqref{defthetaB}.
Part (ii) is a special case of Riemann's singularity theorem (see
e.g. \cite{kempf1}).
The differential at the origin of the separable isogeny $V: J_1 \ra J$ is
an isomorphism $dV: \T_0J_1 \map{\sim} \T_0J$, which identifies
with the Hasse-Witt map $F: H^1(X_1,\cO_{X_1}) \ra H^1(X,\cO_{X})$.
Given $g \not=0$, it
will be enough to compute the tangent space to the
divisor $V^* (T^*_g \Theta)$ at the origin.
Let $\{ Y_g \}_{g \in G}$ be the canonical Theta basis of
$H^0(J_1, \cO(\Theta_1 + T^*_g \Theta_1))$ (Apply Lemma \ref{thetabasis}(ii)
to $L_1 = \cO(\Theta_1 + T^*_g \Theta_1) = \cO(2 T_h^*\Theta_1)$ with
$2h=g$).
Then, by the isogeny formula \eqref{isogenyfor2}, we have
(up to a scalar) $V^* (T^*_g \Theta) = \sum_{g \in G} Y_g$.
Let $\phi_v : \Spec(R) \ra J_1$ be a tangent vector to $J_1$ at the
origin. Then we compute, using $v.Y_g = e^L(v,g) Y_g$, $\forall
v \in \g$
$$ \phi_v^*(\sum_{g \in G} Y_g) = \epsilon e^L(v,g)
\sum_{h \in G/ \langle g \rangle} Y_h (0) $$
Form this we deduce that $\Theta_1$ is singular at $g$ if and only if
$\sum_{h \in G/ \langle g \rangle} Y_h (0) = 0$ and,
assuming $g$ smooth, the equation of $\T_g \Theta_1$. The second
description of $\T_g \Theta$ is a consequence of \cite{kempf1}.
\end{proof}

\begin{que}
Are there other principally polarized abelian varieties $(A,\Theta)$
than Jacobians which have property (i) of Lemma \ref{thetainf}?
\end{que}
\begin{definition}
We say that $X$ has {\em no vanishing theta-null} if $X$ is ordinary
and if all theta characteristics $\kappa$ different from $B$
satisfy $h^0(X,\kappa)=1$, or equivalently (Lemma \ref{thetainf}(ii))
all non-zero $2$-torsion points are smooth points of $\Theta$.
\end{definition}

In the next section we will see (Proposition \ref{propfactorization}(1))
that a generic curve has no vanishing theta-null.

\subsection{Isogeny formulae}

Given an isogeny $f: X \ra Y$ and a line bundle $L$ on $Y$, the
isogeny formula gives the linear map $f^*: H^0(Y,L) \ra H^0(X,f^*L)$
in terms of the canonical Theta bases. Although originally proved
for separable isogenies and line bundles of separable type, we
can extend the isogeny formula to more general line bundles.
We assume $p=2$ and we present (without proof) the three cases
needed in this paper.

\begin{enumerate}
\item
The separable isogeny $V: J_1 \ra J$ with kernel $G$. Let
$\{ X_g \}_{g \in G}$  be the Theta basis of $H^0(J,2\Theta)$ and
$\{ Y^{(2)}_u \}_{u \in G_2}$ be the Theta basis of $H^0(J_1,4\Theta_1)$
(Remark \ref{thetabasispowers}). Then we have
\begin{equation} \label{isogenyV}
\forall g \in G \qquad V^* X_g =
\sum_{u \in G_2 \atop{2u=g}} Y^{(2)}_u.
\end{equation}
\item
Isogeny $V$ as in 1, with $V^*: H^0(J,T^*_g \Theta) \ra
H^0(J_1, \cO(\Theta_1 + T^*_g \Theta_1))$, for $g \in G$.
Then
\begin{equation} \label{isogenyfor2}
V^* X_0 = \sum_{g \in G} Y_g.
\end{equation}
\item
The inseparable isogeny $\xi: J_1 \times J_1 \lra J_1 \times J_1$
defined in \eqref{morxi} with kernel $J_1[2] = G \times \hat{G}$.
Let $A$ be the quotient $J_1 \times J_1/\hat{G}$. Then $\xi$ factorizes
through a separable isogeny (with kernel $G$) $\bar{\xi}: A \ra
J_1 \times J_1$ and we identify $H^0(A,\bar{\xi}^*(2\Theta_1 \boxtimes
2\Theta_1))$ with the $\hat{G}$-invariant subspace of
$H^0(J_1 \times J_1,4\Theta_1 \boxtimes 4\Theta_1)$. A canonical
basis of the latter space is given by the tensors
$\{ Y_u^{(2)} \boxtimes Y_v^{(2)} \}$ with $u,v \in G_2$ such that
$u+v \in G_1 \subset G_2$. The isogeny formula, applied to $\bar{\xi}$,
gives
\begin{equation} \label{isogenyxi}
\forall g,h \in G \qquad \xi^*(Y_h \boxtimes Y_{h+g}) =
\sum_{u,v \in G_2 \atop{\xi(u,v) = (h,h+g)}}
Y^{(2)}_u \boxtimes Y^{(2)}_v.
\end{equation}
\end{enumerate}

\section{Extending Frobenius  to $|2\Theta|$ in characteristic $2$}

We consider a principally polarized Jacobian $(J,\Theta)$ of an
ordinary curve $X$, with $\Theta$ the symmetric Theta
divisor defined by $B$ (section \ref{thetadiv}), and
the morphism $\varphi$ induced by
the linear system $|2\Theta|$ on $J$ (resp. $\varphi_1$ on $J_1$),
which we identify with the Kummer morphism $\kum_J$
(Corollary \ref{kummer}).

\subsection{Factorization}

\begin{prop} \label{propfactorization}
With the notation as above, we have
\begin{enumerate}
\item There exists a non-empty open set of the moduli space of
genus $g$ curves parametrizing curves $X$ with no vanishing
theta-null.

\item Suppose $X$ has no vanishing theta-null. Then the isogeny $V$
can be ``extended''
to a rational map $\tilde{V}$
such that $\varphi \circ V = \tilde{V} \circ \varphi_1$,i.e., the
diagram
\begin{equation} \label{factfrob}
\begin{CD}
 J_1 @>V>> J \\
 @VV\varphi_1 V @VV\varphi V \\
 |2\Theta_1|^* @>\tilde{V}>> |2\Theta|^*
\end{CD}
\end{equation}
is commutative.

\item In terms of the canonical Theta bases of $|2\Theta|$ and
$|2\Theta_1|$ the equations of $\tilde{V}$ are given by $2^g$
quadrics
$$ \tilde{V}: |2\Theta_1|^* \lra |2\Theta|^* \qquad
x:=(x_0:\cdots:x_g:\cdots) \lms (\lambda_0 P_0(x):\cdots:\lambda_g P_g(x):
\cdots)$$
where
\begin{itemize}
\item[(i)] the constants $\{ \lambda_g \}_{g \in G}$, $\lambda_g \in
k$,  satisfy $\forall g \in G, g \not= 0$
$$ \lambda_g = 0  \qquad \iff \qquad  g \
\text{is a singular point of} \  \Theta $$
In particular, if $X$ has no vanishing theta-null, the $\lambda_g$'s
are all non-zero.
\item[(ii)] the polynomials $P_g$ are given by
\begin{equation} \label{polyfact}
  P_g(x) =  \sum_{h \in G/\langle g \rangle} x_{g+h}x_h.
\end{equation}
\end{itemize}
\end{enumerate}
\end{prop}

\begin{proof}
In order to show commutativity of the diagram \eqref{factfrob}
it suffices to show that the image of the injection
$V^*: H^0(J,2\Theta) \hookrightarrow H^0(J_1,4\Theta_1)$
is contained in the image of the multiplication map
\begin{equation} \label{multmap}
\sym^2 H^0(J_1,2\Theta_1) \lra H^0(J_1,4\Theta_1).
\end{equation}

Let  $\{ X_g \}_{g \in G}$ and $\{ Y_g \}_{g \in G}$
be the canonical Theta bases
of $H^0(J,2\Theta)$ and $H^0(J_1,2\Theta_1)$ and consider the
pull-back of equality \eqref{addition}, written for $A=J_1$,
by the morphism $\psi_g: J_1 \ra J_1 \times J_1$, $\psi_g(x) = (x+g,x)$,
\begin{equation} \label{pbadd}
\psi_g^* \xi^*(\theta_1 \boxtimes \theta_1) = \theta_1(g) V^*X_g, \qquad
\psi_g^*(\sum_{h \in G} Y_h \boxtimes Y_h) =
\sum_{h \in G} Y_{h+g} Y_h.
\end{equation}
If $g =0$, we get
\begin{equation} \label{pffeq1}
\sum_{h \in G} Y_h^2 = \theta_1(0) V^*X_0
\end{equation}
with $\theta_1(0) \not= 0$, and define $\lambda_0 = \theta_1(0)$.
If $g \not= 0$,
we see that both members of \eqref{pbadd} are zero. In order to
get a meaningful statement we restrict to a first infinitesimal
neighbourhood of $\psi_g(J_1) \subset J_1 \times J_1$. The
notation is as in section \ref{thetadiv}. We pull-back the
morphisms $\xi$ and $\psi_g$ to $\Spec(R)$ replacing $g$ by
$(1+\epsilon v).g$ and keeping the
same notation for the object over $k$ and its
pull-back to $R$. Pulling back equality \eqref{addition} by
$\psi_g$, for $g \not= 0$, we get the following two
elements in $H^0(J_1,4\Theta_1) \otimes R$
$$ \psi_g^* \xi^* \theta_1 \boxtimes \theta_1 = \epsilon e^L(v,g) \lambda_g
V^*X_g$$
where $\lambda_g$ is a scalar which vanishes if and only if $\Theta_1$ is
singular at the point $g$ (Lemma \ref{thetainf}), and
\begin{eqnarray*}
\psi_g^*(\sum_{h \in G} Y_h \boxtimes Y_h) & = & \sum_{h \in G}
\left[ (1+\epsilon v). Y_{g+h} \right] Y_h \\
 & = & \sum_{h \in G} Y_{g+h} Y_h + \epsilon \sum_{h \in G}
(v.Y_{g+h}) Y_h \\
 & = & \epsilon \sum_{h \in G} e^L(v,g+h)Y_{g+h}Y_h \\
 & = & \epsilon e^L(v,g) \sum_{h \in G/\langle g \rangle} Y_{g+h} Y_h
\end{eqnarray*}
where $h$ runs over a set of representatives of $G/\langle g \rangle$.
Since these elements are equal $\forall v \in \g$ and
there exists $v \in \g$ such that $e^L(v,g) \not= 0$, we obtain
(up to a multiplicative non-zero scalar)
\begin{equation} \label{pffeq2}
\sum_{h \in G/\langle g \rangle} Y_{g+h} Y_h =
\lambda_g V^*X_g.
\end{equation}
In order to complete the proof it suffices to show that for
a general curve $X$ the
$\lambda_g$'s are non-zero, for all $g \not= 0$. Assume that
the contrary holds. Then the equality in $\sym^2 H^0(J_1,2\Theta_1)$
(which we leave as an exercise)
$$ \left( \sum_{h \in G_g} Y_h \right)
\left( \sum_{h \notin G_g} Y_h \right) =
\sum_{h \notin G_g} P_h $$
where $G_g$ is any index $2$ subgroup of $G$ not containing $g$
and the $P_h$ are the quadrics \eqref{polyfact}, shows that
either $\sum_{h \in G_g} Y_h$ or
$\sum_{h \notin G_g} Y_h$ is zero, since the RHS is mapped to zero
in $H^0(J_1,4\Theta_1)$. But this is impossible, since $\{ Y_h \}$ is a
basis. Hence for any curve $X$ there exists $g \not= 0$ such that
$\lambda_g \not= 0$ and we conclude by a monodromy
argument found in \cite{Eke}.

\end{proof}

\begin{rem}
We notice that $\tilde{V}$ is uniquely defined only up to a degree $2$
equation of the image $\varphi_1(J_1) \subset |2\Theta_1|^*$.
We will show uniqueness of $\tilde{V}$ (Proposition \ref{equationkummer})
for an ordinary genus $2$ curve.
We expect that there are no quadrics containing $\varphi_1(J_1)$ for
a curve of genus $g \geq 2$ with no vanishing theta-null.
\end{rem}

Equalities \eqref{pffeq1} and \eqref{pffeq2} can be used to
define the vector $(\lambda_g)_{g \in G}$ up to a scalar. In the
next proposition we give a more direct definition in terms of
theta-constants. Let $\{ Y_u^{(2)} \}_{u \in G_2}$ be the Theta
basis of $H^0(J_1,4\Theta_1)$ and we denote by $Y_u^{(2)}(0) \in k$
the value of $Y_u^{(2)}$ at the origin (after having chosen an
isomorphism $\cO(4\Theta_1)_0 \map{\sim} k$).

\begin{prop}
With the notation as above, we have
\begin{equation} \label{exprlambda}
\forall g \in G \qquad \lambda_g =
\sum_{u \in S_g} Y_u^{(2)}(0)
\end{equation}
where $S_g$ is a set of representatives of
$\{u \in G_2 \ : \  2u=g\}/\langle g \rangle$.
\end{prop}

\begin{proof}
Let $i$ be the inclusion in the first factor $i: J_1 \hookrightarrow
J_1 \times J_1$. A standard computation modelled on \cite{kempf} Proposition
6, which involves \eqref{isogenyxi}, gives the equality in
$\sym^2 H^0(J_1,4\Theta_1)$
$$ i^* \xi^* \left( \sum_{h \in G/\langle g \rangle}
Y_h \boxtimes Y_{h+g} \right) = \left(\sum_{u \in S_g}
Y_u^{(2)}(0) \right) \cdot  \left( \sum_{u \in G_2 \atop{2u=g}}
Y_u^{(2)} \right) $$
We observe that the composite $\xi \circ i$ is the diagonal map,
hence the LHS can be rewritten as $\sum Y_h Y_{h+g}$. Applying
the isogeny formula \eqref{isogenyV}, we get $\forall g \in G$
$$ \sum_{h \in G/\langle g \rangle}
Y_h Y_{h+g} = \left(\sum_{u \in S_g}
Y_u^{(2)}(0) \right) \cdot V^* X_g.$$
Comparing with \eqref{pffeq1} and \eqref{pffeq2} we can conclude.
Finally we observe that the expression \eqref{exprlambda}
is well-defined since $\forall u$ such that
$2u=g$ we have $Y_u^{(2)}(0) = Y_{-u}^{(2)}(0) = Y_{u+g}^{(2)}(0)$
($Y_0^{(2)}$ is symmetric).
\end{proof}

\section{Kummer's quartic surface in characteristic $2$}

As an application of Proposition \ref{propfactorization}
we shall deduce the
equation of the image $\varphi(J) \subset |2\Theta|^*$
for a Jacobian of an ordinary genus $2$ curve $X$. We observe that by
Clifford's theorem $h^0(X_1,Bg) \leq 1$, $\forall g \in G$, hence $X$
has no vanishing theta-null and we can apply Proposition
\ref{propfactorization}(2).
Somewhat surprisingly, invariance under the Theta group $G(L)$ and
the ``Frobenius'' map $\tilde{V}$ are sufficient to determine
the equation.
We fix an isomorphism $G \cong (\ZZ/2\ZZ)^2$.

\begin{prop} \label{equationkummer}
In the canonical Theta basis $\{X_g\}_{g \in G}$ of $H^0(J,2\Theta)$
the equation of the Kummer surface is (up to a non-zero scalar)
\begin{equation} \label{kqs}
\lambda^2_{10}(x_{00}^2x_{10}^2 + x_{01}^2 x_{11}^2) +
\lambda^2_{01}(x_{00}^2x_{01}^2 + x_{10}^2 x_{11}^2) +
\lambda^2_{11}(x_{00}^2x_{11}^2 + x_{01}^2 x_{10}^2) +
\frac{\lambda_{10} \lambda_{01} \lambda_{11}}{\lambda_{00}}
x_{00}x_{10}x_{01}x_{11}
\end{equation}
In particular, $\varphi(J)$ is not contained in a quadric and
therefore the map $\tilde{V}$ is uniquely defined.
\end{prop}

\begin{proof}
First we observe that the image $\varphi(J) \subset |2\Theta|$ is
invariant under the Theta group $G(L)$. Hence the equation
of $\varphi(J)$ has to be $G(L)$-invariant. Since $(2\Theta)^2 = 8$
and $\deg \varphi = 2$ ($\varphi$ is separable), we have $\deg
\varphi(J) = 4$. The unique $G(L)$-invariant quadric
is $x_{00}^2 + x_{01}^2 + x_{10}^2 + x_{11}^2 =
(x_{00} + x_{01} + x_{10} + x_{11})^2 $. Since $\varphi(J)$ is
non-degenerate, we see that it is not contained in a
quadric. A straightforward computation shows that a basis
of $G(L)$-invariant quartics is given by the $5$
polynomials
$$ S:= x^4_{00} + x^4_{01} + x^4_{10} + x^4_{11} \qquad
R:= x_{00}x_{01}x_{10}x_{11}$$
$$ Q_{10} := x^2_{00}x^2_{10} + x^2_{01}x^2_{11} \qquad
Q_{01} := x^2_{00}x^2_{01} + x^2_{10}x^2_{11} \qquad
Q_{11} := x^2_{00}x^2_{11} + x^2_{01}x^2_{10} $$
Let us denote by
$$ F= \alpha_{10} Q_{10} + \alpha_{01} Q_{01} + \alpha_{11} Q_{11}
+ \beta R + \gamma S $$
the equation of $\varphi(J)$, where $\alpha_{10}, \alpha_{01},
\alpha_{11}, \beta, \gamma \in k$ need to be determined.
First $\varphi(e) = (1:0:0:0)$, where $e \in J$ is the origin,
implies that $\gamma = 0$. Since the image is reduced, $\beta \not= 0$.
Next, the equation of $\varphi_1(J_1)
\subset |2\Theta_1|$ is given in a Theta basis by
$$ F_1 = \alpha^2_{10} Q_{10} + \alpha^2_{01} Q_{01} + \alpha^2_{11}
Q_{11} + \beta^2 R $$
Since $\varphi_1(J_1)$ is not contained in a quadric, the rational
map $\tilde{V}$ \eqref{factfrob} is uniquely defined. Since
$\tilde{V}(\varphi_1(J_1)) = \varphi(J)$, there exists a
$G(L)$-invariant quartic $A$, such that (up to a non-zero scalar)
\begin{equation} \label{invquar}
\tilde{V}^*(F) = F_1 \cdot A
\end{equation}
We shall use the Theta coordinates $\{ x_g \}_{g \in G}$ on both
spaces. In order to determine the equation of $A$, we restrict
equality \eqref{invquar} to the hyperplane $H: \ \sum_{g \in G}
x_g = 0$ and write the expressions as degree $8$ polynomials
in $x_{01},x_{10},x_{11}$. A straightforward computation, which
we omit, leads to
\begin{eqnarray*}
\tilde{V}^*(F)_{|H} = \lambda^2_{01} \lambda^2_{10} \lambda^2_{11}
&&(x_{01} + x_{11})^2 (x_{01} + x_{10})^2 (x_{11} + x_{10})^2 \\
&&\left[ x_{01}^2 (\frac{\alpha_{10}}{\lambda^2_{10}} +
\frac{\alpha_{11}}{\lambda^2_{11}})
+  x_{10}^2 (\frac{\alpha_{01}}{\lambda^2_{01}} +
\frac{\alpha_{11}}{\lambda^2_{11}})
+  x_{11}^2 (\frac{\alpha_{10}}{\lambda^2_{10}} +
\frac{\alpha_{01}}{\lambda^2_{01}}) \right]
\end{eqnarray*}
Suppose that $A_{|H} \not= 0$, i.e. $\tilde{V}^*(F)_{|H} \not= 0$.
Then at least one of the factors $x_{01} + x_{10}, x_{01} + x_{11},
x_{11} + x_{10}$ has to divide ${F_1}_{|H}$. Again elementary
computation shows that this can only happen when $\beta = 0$,
which is impossible. Hence $A$ is a multiple of $H$ and, since
$A$ is also $G(L)$-invariant, we get $A= H^4 = \sum_{g \in G} x_g^4$.
Moreover $\tilde{V}^*(F)_{|H} = 0$ implies that there exists
a $\mu \in k^*$ such that $\alpha_g = \mu \lambda_g^2$ for all
$g \in G^*$.
\bigskip

It remains to determine the constant $\beta$. We compute the
degree $8$ polynomial $\tilde{V}^*(F)$ (again we omit details)
\begin{align*}
\mu\lambda^2_{00} H^4 (\lambda^4_{10} Q_{10} + \lambda^4_{01}
Q_{01} + \lambda^4_{11} Q_{11}) \\
+\mu \lambda^2_{01} \lambda^2_{10} \lambda^2_{11} H^2
(x^2_{00}x^2_{01}x^2_{10} + x^2_{00}x^2_{01}x^2_{11} +
x^2_{00}x^2_{10}x^2_{11} + x^2_{01}x^2_{10}x^2_{11}) \\
+ \beta \lambda_{00} \lambda_{01} \lambda_{10} \lambda_{11}
H^2 (RH^2 + x^2_{00}x^2_{01}x^2_{10} + x^2_{00}x^2_{01}x^2_{11} +
x^2_{00}x^2_{10}x^2_{11} + x^2_{01}x^2_{10}x^2_{11})
\end{align*}
This expression being equal to $F_1 \cdot H^4$, we get the equation
mentioned in the proposition.
\end{proof}

\begin{rem}
If the characteristic is different from $2$, Kummer's quartic surface is a
much studied object (see e.g.\cite{hud}). Among many other
results, let us just mention that the coefficients of the
quartic equation (in a Theta basis) satisfy a cubic equation
and are themselves polynomials of degree $12$ in the
theta-constants. Since we could not find any treatment
of the characteristic $2$ case in the literature, we decided
to include it in our paper.
\end{rem}

\section{The moduli space $\mod0$ of rank $2$ vector bundles}

We assume $p=2$.
Let $\mod0$ denote the moduli space of rank $2$ semi-stable
vector bundles over $X$ with trivial determinant.
As was observed in \cite{MehRam}, the theta divisor
$$\tilde{\Theta} =\{[E]\in \mod0 \ | \ h^0(E\otimes B)\not=0\}$$
is Cartier (and not only $\QQ$-Cartier) and is ample by GIT. We
denote by $\LL_0$ the line bundle $\cO_{\mod0}(\tilde{\Theta})$.
By \cite{Ray} there exists a regular morphism
$$D:\ \mod0 \lra \PP H^0(J,2\Theta)=|2\Theta|$$
which maps the class of the semi-stable bundle $E$ to the divisor $D(E)$
$$\supp D(E)=\{ L \in J \ |\ H^0(E\otimes B\otimes L)\not= 0\}.$$
As in the complex case \cite{NarRam} one has
\begin{prop}[V. Balaji] \label{isomxp3}
If $g=2$, the morphism $D: \mod0 \lra \pp^3 = |2\Theta|$ is
an isomorphism.
\end{prop}
\begin{proof}
(sketch) We proceed in two steps. First we consider a flat family
of curves $\XX \lra T = \Spec(A)$, where $A$ is a discrete
valuation ring such that its residue field at the closed point
$0$ is $k$, its field of fraction $K$ is of characteristic zero,
and $\XX_0 = X$. We consider the moduli scheme $\M \lra T$ of
semi-stable rank $2$ vector bundles of trivial determinant
over the family $\XX \ra T$. Let $\M_0$ be the fibre of
$\M \lra T$ over $0$. Then, by GIT over arbitrary base \cite{ses},
we have a canonical bijective morphism $i: \mod0 \lra
\M_0$. Moreover on the open set
of stable points $i$ is an isomorphism since the action
of the projective group on the Quot scheme is free.
We conclude that $i: \mod0 \lra \M_0$ is an
isomorphism by Zariski's Main Theorem.

\bigskip

Secondly, we extend the morphism $D$ \cite{NarRam} to the
family $\M \lra T$,i.e., we construct a morphism over the base $T$
$$ \DD: \M \lra \pp^3_T $$
such that $\DD_0 = D$. In order to show the proposition, it will
be enough (again by Zariski's Main Theorem) to show that
$\DD$ is birational and bijective. We consider the fibre
$\DD_\xi: \M_\xi \lra \pp^3_\xi$ over the generic point
$\xi \in T$. Working over an algebraic closure of $K$ (of characteristic
$0$), we see \cite{NarRam} that $D_\xi$ is an isomorphism.
Hence $\DD$ is birational. It remains to show that $\DD_0 = D$
is bijective. Surjectivity is obvious.
Let $H$ be the hyperplane $|2\Theta|$ of divisors
passing through $\cO$. The inverse image of $H$ is $\Theta$ showing
$D^*\cO(1)\cong \LL_0$. It follows that $D$ is finite
because $D^*\cO(1)$ is ample.
Let $\LL = D^* \cO_{\pp^3_T}(1)$ be the relatively ample
line bundle over $\M \lra T$. Consider the canonical inclusion
of $\cO_{\pp^3_T}$ in $\DD_* \cO_{\M}$ with cokernel $\mathcal{Q}$.
$$ 0 \lra \cO_{\pp^3_T} \lra  \DD_* \cO_{\M} \lra \mathcal{Q} \lra 0 $$
We twist by $\cO_{\pp^3_T}(n)$. If $n$ is large enough, we have
$\forall t \in T$ and $\forall i > 0$,
$h^i(\M_t,\LL_t^n) = 0$. Hence, since $\pp^3_T \lra T$ is flat, we
see that
$$ h^0(\pp^3_t,\cO_{\pp^3_t}(n)) - h^0(\M_t,\LL_t^n) =
\chi(\pp^3_t,\cO_{\pp^3_t}(n)) - \chi(\M_t,\LL_t^n) $$
is constant. If $t = \xi$, this number is zero \cite{NarRam}.
Hence $h^0(\pp^3_0,\mathcal{Q}_0(n)) = 0$ and $Q_0 = 0$. So
$\cO_{\pp^3_0} = D_* \cO_{\mod0}$ and $D$ is injective.

\end{proof}

\begin{rem}
As we were told  by C.S. Seshadri, the first part of the proof is
completely worked out in the PhD thesis of Venkata Balaji (in preparation).
A direct proof of this isomorphism along the lines of the
original paper \cite{NarRam} was obtained by M.S. Narasimhan.
\end{rem}

\section{Frobenius action on $\mod0$ for an ordinary genus $2$ curve}

The goal of this section is to describe the Frobenius map
(more precisely, its separable part, the Verschiebung)
$$V: \MM_{X_1} \lra \mod0 , \qquad E \lms F^* E. $$
We consider the morphism $\psi: J \ra \mod0$, $L \mapsto
[L \oplus L^{-1}]$ which, when composed with $D$, equals
the Kummer morphism $\kum_J$. Because of Proposition \ref{isomxp3}
and Corollary \ref{kummer} we can identify $\psi$ and $\varphi$.
Since $H^0(\mod0,\LL_0) = H^0(J,2\Theta)$ and since $\tilde{V}$ is
uniquely defined (Proposition \ref{equationkummer}), we can
identify (via $D$) the Verschiebung $V: \MM_{X_1} \lra \mod0$ with
the rational map $\tilde{V}: |2\Theta_1|^* \lra |2\Theta|^*$ given
by the equations \eqref{polyfact}.
We gather our results in the following proposition.
\begin{prop} \label{computations}
Let $X$ be an ordinary genus $2$ curve.
\begin{enumerate}
\item
The semi-stable boundary of $\mod0$ (resp. $\MM_{X_1}$) is
isomorphic to Kummer's quartic surface $\varphi(J)$ (resp.
$\varphi_1(J_1)$) whose equation is given in \eqref{kqs}.
In particular, $V$ maps $\varphi_1(J_1)$ onto $\varphi(J)$.
\item
There exists a unique stable bundle $E_{BAD} \in \MM_X$, which
is destabilized by the Frobenius map, i.e. $F^* E_{BAD}$ is
not semi-stable. We have $E_{BAD} = {F}_* B^{-1}$ and its projective
coordinates are $(1:1:1:1)$.
\item
The set of bundles $\{ [E] \in \MM_{X_1} \ | \ \Hom(E, F_{*}B)
\not= 0 \}$ is the hyperplane $H_1: \sum_{g \in G} x_g = 0$.
In particular, $E_{BAD} \in H_1$. The restriction of $V$ to $H_1$
contracts $H_1$ to the conic $\varphi(J) \cap H$, where $H$ is the
hyperplane $x_0 = 0$ in $|2\Theta|$. In particular, any stable
bundle $E \in H_1$ is mapped  into the semi-stable boundary of
$\mod0$.
\item
The fiber of $V$ over a point $[E] \in \MM_X$ is
\begin{itemize}
\item
a non-degenerate $G$-orbit of a point $[E_1] \in \MM_{X_1}$ ($4$
distinct points), if $[E] \notin H$
\item
empty, if $[E] \in H \setminus (H \cap \varphi(J))$
\item
a projective line passing through $E_{BAD}$,
if $[E] \in H \cap \varphi(J)$
\end{itemize}
In particular, $V$ is not surjective and the separable degree of $V$
is $4$.
\end{enumerate}
\end{prop}

\begin{proof}
1. This follows immediately from \cite{NarRam} and Proposition
\ref{equationkummer}.

2. The base locus of $\tilde{V}$ is given by the intersection
$\bigcap_{g \in G} \{ P_g = 0 \}$, which turns  out (after some
elementary computations) to be a unique point with
projective coordinates $(1:1:1:1)$. In terms of vector bundles
this point, denoted $E_{BAD}$, corresponds to the direct image
$F_{*} B^{-1}$. Indeed, $F_{*} B^{-1}$ is stable: a nonzero map
$L \ra F_{*} B^{-1}$ is equivalent, by adjunction, to a non-zero
map $F^* L \ra B^{-1}$, hence $2 \deg L \leq -1$. Moreover,
since we have a canonical nonzero map $F^*F_{*} B^{-1} \ra B^{-1}$,
this bundle is destabilized. Uniqueness can also be proved
without the use of the equations: consider any stable bundle $E
\in \MM_X$, which is destabilized by $V$. By \cite{ls} Corollary 2.6
there exists a theta-characteristic $A$ which appears as
a quotient $F^*E \ra A^{-1}$. By adjunction this quotient
map induces a map $E \ra F_*A^{-1}$. Since the two bundles
are stable with the same slope, we deduce that they have to
be isomorphic. The determinant of $F_* A^{-1}$ being trivial,
we get $A = B$ and we are done.

3. First we observe that the hyperplane $H_1$ is mapped into the
hyperplane $H$. After fixing an isomorphism $G \cong (\ZZ/2\ZZ)^2$,
we straightforwardly compute that a point in the image $V(H_1)$ satisfies
the equations
\begin{equation} \label{inthkqs}
\lambda_{10} x_{01} x_{11} + \lambda_{01} x_{10} x_{11} +
\lambda_{11} x_{01} x_{10} = 0
\end{equation}
which is precisely (after squaring) the equation of Kummer's
quartic surface \eqref{kqs} restricted to the hyperplane $H$. Moreover
$V^*(H) = H_1^2$ and by adjunction $\Hom(E,F_{*}B) = \Hom(F^* E,B)$
which proves the first assertion.

4. Given a point $[E] \in \MM_X$ with projective coordinates
$(a_{00}:a_{01}:a_{10}:a_{11})$, we have to solve the
system of quadratic equations
\begin{equation} \label{syseq}
P_g = \frac{a_g}{\lambda_g} \qquad \forall g \in G \cong
(\ZZ/2\ZZ)^2
\end{equation}
where the quadrics $P_g$'s are defined in \eqref{polyfact}.
We write $b_g = \frac{a_g}{\lambda_g}$. Adding any two $P_g$'s
with $g \not= 0$, we find
\begin{equation} \label{eq1}
P_{01} + P_{10} = x_{00}x_{01} + x_{10}x_{11} + x_{00}x_{10} +
x_{01}x_{11} = (x_{00} + x_{11})(x_{01} + x_{10}) = b_{01} +
b_{10}
\end{equation}
and similarly
\begin{equation} \label{eq2}
(x_{00} + x_{01})(x_{11} + x_{10}) = b_{11} + b_{10} \qquad
(x_{00} + x_{10})(x_{01} + x_{11}) = b_{01} + b_{11}
\end{equation}
and
\begin{equation} \label{eq3}
x_{00} + x_{01} + x_{10} + x_{11} = c
\end{equation}
with $c^2 = b_{00}$. We let $\alpha = x_{00} + x_{11},
\beta = x_{01} + x_{10}, \gamma = x_{00} + x_{01},
\delta = x_{11} + x_{10}$. Then the equations \eqref{eq1},
\eqref{eq2}, \eqref{eq3} imply that $\alpha, \beta$ (resp.
$\gamma,\delta$) are the roots of the polynomial
\begin{equation} \label{poly}
t^2 + c t + (b_{01} + b_{10}) \qquad \text{resp.}\ \
t^2 + c t + (b_{11} + b_{10})
\end{equation}
Substituing $x_{01},x_{11},x_{01}$ with $\gamma + x_{00},
\alpha + x_{00}$ and $\beta + \gamma + x_{00}$ respectively in
the equation $P_{01} = b_{01}$, we find
\begin{equation} \label{eq4}
cx_{00} = b_{10} + \gamma \alpha
\end{equation}
Assuming $c\not= 0$,i.e. $[E] \notin H$, we see that
the fiber $V^{-1}([E])$ consists of the point $[E_1] \in
\MM_{X_1}$ with projective coordinates $(x_{00}:x_{01}:x_{10}:x_{11})$
given by
$$ x_{00} = \frac{1}{c}(b_{10} + \gamma \alpha) \qquad
x_{01} = \frac{1}{c}(b_{10} + \gamma \beta) \qquad
x_{10} = \frac{1}{c}(b_{10} + \beta \delta) \qquad
x_{11} = \frac{1}{c}(b_{10} + \delta \alpha)$$
plus the other $3$ points obtained by switching $\alpha$ and
$\beta$ as well as $\gamma$ and $\delta$. Since $c \not= 0$,
one easily sees that these $4$ points are distinct.

Assume $c=0$, i.e. $[E] \in H$. We get $\alpha = \beta, \gamma =
\delta, \alpha^2 = b_{01} + b_{10}, \gamma^2 = b_{11} + b_{10}$.
Hence $(b_{10} + \gamma \alpha)^2 = b_{10} b_{11} + b_{01} b_{10}
+ b_{10} b_{11}$. Assume that $b_{10} b_{11} + b_{01} b_{10}
+ b_{10} b_{11} \not= 0$,i.e. $[E] \in
H \setminus (H \cap \varphi(J))$. Because of \eqref{eq4} the
system \eqref{syseq} does not have a solution. If
$ b_{10} b_{11} + b_{01} b_{10} + b_{10} b_{11} = 0$, i.e.
$[E] \in H \cap \varphi(J)$, it is easy to check that
$V^{-1}([E])$ is the $\PP^1$, intersection of the two hyperplanes
$$ H_1: x_{00} + x_{10} + x_{01} + x_{11} = 0 \qquad
H_{\alpha,\gamma}: (\alpha + \gamma)x_{00} + \alpha x_{01}
+ \gamma x_{11} = 0 $$
\end{proof}

\begin{rem}
Let $L_1 \in JX_1$ and $M = F^*L_1 \in JX$. We assume that
$M^2 \not= \cO_X$. Then obviously $F^*[L_1 \oplus L_1^{-1}]
= [M \oplus M^{-1}]$ and, as one easily checks, if
$h^0(M \otimes B) = 0 $, the three
isomorphism classes contained in  $[L_1 \oplus L_1^{-1}]$, i.e. the
decomposable bundle $L_1 \oplus L_1^{-1}$, the two non-split
extensions of $L_1$ by $L_1^{-1}$ and of $L_1^{-1}$ by $L_1$
(which are interchanged by the hyperelliptic involution $i$),
are mapped to the corresponding isomorphism class in
$[M \oplus M^{-1}]$. On the other hand, if $h^0(M \otimes B) >0$,
all three isomorphism classes are mapped to $M \oplus M^{-1}$.
Moreover, in that case,  by Proposition
\ref{computations} (4), there exist stable bundles $E$ such
that $[F^*E] = [M \oplus M^{-1}]$. Since $i$ commutes
with $F$ and since any stable bundle on $X_1$ is $i$-invariant,
we have $F^*E \cong M \oplus M^{-1}$. This shows that the
non-split extension of $M$ by $M^{-1}$ is not of the
form $F^*E$ for some semi-stable bundle $E$.
\end{rem}

\begin{rem}
It would be interesting to have a coordinate-free
proof of the following fact, which follows from Proposition
\ref{computations}(3):
If $E_1$ is stable over $X_1$ with
$\Hom(E_1,F_{*}B) \not= 0$, then $F^* E_1$ is unstable.
\end{rem}

Let $\NN_X$ (resp. $\NN_{X_1}$) denote the moduli space of
semi-stable rank $2$ vector bundles of degree $0$ over $X$
(resp. $X_1$). As a corollary of Proposition \ref{computations}
we have
\begin{prop} \label{surjectivity}
For any ordinary curve $X$, the rational map
$V: \NN_{X_1} \lra \NN_{X}$ given by $E \lms F^* E$
is surjective.
\end{prop}

\begin{proof}
It will be enough to show that any point $[E] \in \MM_X \subset \NN_X$
is in the image, since, twisting by a degree zero line bundle, we
can always assume the determinant to be trivial. For any non-zero
$g \in G$, we choose a $h \in G_2$ such that $2h = g$ and we
denote by $\MM_{X_1}(g)$ the moduli of rank $2$ vector bundles
with fixed determinant equal to $g$. Then we have
a commutative diagram
$$
\begin{CD}
\MM_{X_1}(g) @>V>> \MM_X \\
@VVT_hV  @VVT_gV \\
\MM_{X_1} @>V>> \MM_X
\end{CD}
$$
Now since the vertical maps are isomorphisms, the image of
the first horizontal map contains the complement to the hyperplane
$x_g \not=0$. Since the hyperplanes $\{ x_g = 0 \}_{g \in G}$ have
no common point in $\MM_X$ we get surjectivity.
\end{proof}

\section{The action on the odd degree moduli space $\MM_{X_1}(\Delta)$}

In this section we briefly discuss the action of the Frobenius
map on the moduli space $\MM_{X_1}(\Delta)$ of semi-stable
rank $2$ vector bundles over $X_1$ with fixed determinant
equal to $\De$, with $\deg \De = 1$.
We will study the rational map
$$ V: \MM_{X_1}(\De) \lra \mod0, \qquad E \lms F^* E \otimes \De^{-1} $$
Note that we use the same letter ($\De,B,\ldots$) for
the line bundles over $X_1$ and $X$, which correspond under
the $k$-semi-linear isomorphism $\iota: X_1 \ra X$.

The next proposition holds for any curve of genus $g \geq 2$.

\begin{prop}
The image of $V$ is contained in the Theta divisor $\Theta \subset
\mod0$, i.e.
$$\forall E \in \MM_{X_1}(\Delta) \qquad
h^0(X,F^*E \otimes \De^{-1} \otimes B) >0$$
In particular, $V$ is not dominant.
\end{prop}

\begin{proof}
We write $E$ as an extension of line bundles
$$ 0 \lra L \lra E \lra \De L^{-1} \lra 0 \qquad (\epsilon)$$
for some line bundle $L$ with $\deg L \leq 0$ and $\epsilon
\in \pp H^1(X_1,L^2 \De^{-1})$. The extension class $F^* \epsilon$
of the exact sequence, gotten by pull-back under $F$, i.e.
\begin{equation} \label{eqpr1}
0 \lra L^2 \lra F^* E \lra \De^2 L^{-2} \lra 0 \qquad
(F^* \epsilon)
\end{equation}
is obtained from $\epsilon$ via the linear map
\begin{equation} \label{eqpr2}
 H^1(X_1,L^2 \De^{-1}) \map{F^*} H^1(X,L^4\De^{-2}) =
H^1(X_1,L^2\De^{-1} \otimes F_* \cO_X)
\end{equation}
The last map coincides with the induced map on
cohomology of the canonical exact sequence
\begin{equation} \label{eqpr3}
0 \lra L^2 \De^{-1} \lra L^2 \De^{-1} \otimes F_* \cO_X \lra
L^2 \De^{-1} B \lra 0
\end{equation}
In order to prove that $h^0(X, F^*E \otimes \De^{-1}B) >0$,
we tensorize \eqref{eqpr1} with $\De^{-1} B$. If
$h^0(X,L^2 \De^{-1} B) >0 $, we are done. Therefore we assume
$h^0(X,L^2 \De^{-1} B) = 0$. It follows from \eqref{eqpr2}
and \eqref{eqpr3} that
\begin{equation} \label{eqpr4}
0 \lra H^1(X_1,L^2 \De^{-1}) \map{F^*} H^1(X, L^4 \De^{-2}) \lra
H^1(X_1, L^2 \De^{-1} B) \lra 0
\end{equation}

On the other hand, we see that
$h^0(X, F^*E \otimes \De^{-1} \otimes B) >0$ if and only if
the symmetric coboundary map
$$ H^0(X, \De L^{-2} B) \map{\cup F^* \epsilon}
H^1(X,\De^{-1} L^2 B) = H^0(X,\De L^{-2} B)^*$$
is degenerate. We write $V := H^0(X, \De L^{-2} B)$. It is
well-known that the linear map
\begin{equation} \label{eqpr5}
H^1(X,L^4 \De^{-2}) = H^0(X, \Omega_X L^{-4} \De)^* \map{m^*}
V^* \otimes V^*, \qquad \delta \lms \cup \delta
\end{equation}
is the dual of the multiplication map of global sections
$\sym^2 V \map{m} H^0(X, \Omega_X L^{-4} \De)$.
Let us denote by $\DD_2(V^*)$ the space of divided
powers of $V^* \otimes V^*$,i.e. the subspace of
tensors invariant under the involution $\phi \otimes \psi \lms
\psi \otimes \phi$. We have $\DD_2(V^*) = (\sym^2 V)^*$.
We denote by $F(V)$ the subspace of $\sym^2 V$ generated by the
squares $v^2$ with $v \in V$. Dually, the kernel of the
surjection $ K:= \ker ( \DD_2(V^*) \lra F(V)^* )$
coincides with the space of alternating maps $V \ra V^*$.
The main point,
which will be used later, is that, since, by Riemann-Roch, $\dim V
= \dim H^1(X,\De^{-1} L^2 B) = 1 - 2\deg L$ is odd, any
map in $K$ is degenerate.

By Serre duality we have
$H^1(X_1,L^2 \De^{-1} B) = H^0(X_1,\De L^{-2} B)^*$. Put
$V_1:= H^0(X_1,\De L^{-2} B)$. Then we have a $k$-linear isomorphism
$\varphi: V_1 \map{\sim} F(V)$ defined as follows: we have
$V_1 = \iota^* V = V \otimes_k k$ and we put $\varphi(v \otimes t)
=tv^2 \in F(V)$. We also note that $F(V)^* = F(V^*)$.
Again by Serre duality, we observe that the linear maps \eqref{eqpr4}
and \eqref{eqpr5} coincide, i.e. we have a commutative
diagram
$$
\begin{array}{ccccccccc}
 H^1(X_1, L^2\De^{-1}) & \map{F^*} & H^1(X,L^4\De^{-2})
& & \lra & & H^1(X_1, L^2 \De^{-1} B) & \lra & 0 \\
  & & \downarrow & & & &\downarrow & &  \\
   & & H^0(X, \Omega_X L^{-4} \De)^* & \map{m^*}& \DD_2(V^*) &
\lra & F(V)^* & &
\end{array}
$$
where the vertical maps are $k$-linear isomorphisms. Now
we can conclude as follows: by commutativity, any extension
class $F^* \epsilon$, with $\epsilon \in H^1(X_1,L^2 \De^{-1})$,
is mapped by $m^*$ into $K$. Hence the corresponding coboundary
map is degenerate.
\end{proof}

We come back to an ordinary curve $X$ of genus $2$.
Let us denote by $\Gr \subset \pp \Lambda^2 H^0(J_1,2\Theta_1):= \pp^5$
the Grassmannian of projective lines in $\pp H^0(J_1,2\Theta_1) = \pp^3 =
\MM_{X_1}$. Following \cite{beau} section 3.4 we consider, for a
general point $q \in X$, the morphism
$$\MM_{X_1}(\cO(q)) \map{\L} \Gr, \qquad E \lms \L (E) =
\{ E' \in \MM_{X_1}: \ E' \subset E \}. $$
We choose an isomorphism $G \cong (\ZZ/2\ZZ)^2$ and
consider the Pl\"ucker coordinates on $\Lambda^2 H^0(J_1, 2\Theta_1)$
\begin{eqnarray*}
z_1 = x_{00} \wedge x_{01}, \qquad z_2 = x_{00} \wedge x_{10}, \qquad
z_3 = x_{00} \wedge x_{11}, \\
z_4 = x_{10} \wedge x_{11}, \qquad z_5 = x_{01} \wedge x_{11}, \qquad
z_6 = x_{01} \wedge x_{10}.
\end{eqnarray*}
The equation of the Grassmannian $\Gr$ is $z_1z_4 +z_2z_5 + z_3z_6 = 0$.
and the $G$-invariant subspace $Z \subset \Lambda^2 H^0(J_1,2\Theta_1)$
is given by the $3$ linear equations
$$ Z: \qquad z_1+z_4 = z_2 +z_5 = z_3 + z_6 = 0 $$

\begin{prop}
For a general point $q$,
we have a commutative diagram
$$
\begin{CD}
\MM_{X_1}(\cO(q)) @>V>> \MM_X \\
@VV \L V  @VV \cong V \\
\pp \Lambda^2 H^0(J_1,2\Theta_1) @> \tilde{V} >> \pp H^0(J,2\Theta) =
\pp^3
\end{CD}
$$
where $\tilde{V}$ is the projection with center $\pp Z = \pp^2$.
In terms of canonical coordinates on both spaces, we have
$$ \tilde{V}^*(x_{00}) = 0, \qquad
\tilde{V}^*(x_{10}) = z_2 + z_5, \qquad
\tilde{V}^*(x_{11}) = z_3 + z_6, \qquad
\tilde{V}^*(x_{01}) = z_1 + z_4.
$$
There exists a unique bundle in $\MM_{X_1}(\cO(q))$ which is destabilized
by $F$, namely $F_*(B^{-1}(q))$.
\end{prop}

\begin{proof}
Since the proof is in the same spirit as the proof of Proposition
\ref{computations}, we
just give a sketch. Let $D_q$ be the vector field on $J_1$ associated
to $q$. We also denote by $D_q$ the endomorphism of
$H^0(J_1,4\Theta_1)$ obtained via the canonical $p$-integrable
connection $\nabla$
(Remark \ref{Katz}). We observe that we have a commutative
diagram
\begin{equation} \label{wahlmap}
\begin{CD}
\sym^2 H^0(J_1,2\Theta_1) @>>> \Lambda^2 H^0(J_1,2\Theta_1) \\
@VVmV   @VVW_{D_q}V \\
H^0(J_1,4\Theta_1) @>D_q>> H^0(J_1,4\Theta_1)
\end{CD}
\end{equation}
The first horizontal map is the canonical projection and
$W_{D_q}$ is the Wahl map associated to $D_q$ (\cite{beau} section A.10).
We consider the map $J_1 \ra \MM_{X_1}(\cO(q))$ defined in
\cite{beau} section 3 and compose with $\L$. For $q$ general, the
composite is non-degenerate and the induced (injective) map
on global sections coincides with $W_{D_q}$ \cite{beau}. Using
\eqref{wahlmap} we can now deduce the equations of $\tilde{V}$.
The last assertion can be proved as in Proposition
\ref{computations} (2).
\end{proof}

\section{Frobenius dynamics}

Let $F_a$ be the absolute Frobenius map of $X$.
We write  $F_a^{(n)}$ for the $n$-fold composite $F_a \circ \cdots \circ F_a$.
We will study the set of Frobenius semi-stable bundles, i.e.
$$ \Omega^{Frob} : =
\{ [E] \in \mod0 \ | \ F_a^{(n)*}E \ \text{semi-stable} \ \forall
n \geq 1 \} $$
and the set of bundles coming from representations of the
algebraic fundamental group of $X$,i.e. (see \cite{ls} Satz 1.4)
$$ \Omega^{Rep} : =
\{ [E] \in \mod0 \ | \ \exists n>0 \ \ F_a^{(n)*}E \map{\sim} E \} $$
We obviously have $\Omega^{Rep} \subset \Omega^{Frob}$ and we show

\begin{prop} \label{dynamics}
The set $\Omega^{Frob}$ is Zariski dense in $\mod0 = \pp^3$.
\end{prop}

\begin{proof}
Since $F_a = \iota \circ F$, the action of $F_a$ on $\mod0 = \pp^3$
factorizes as
$$ F^*_a : \MM_X \map{i^*} \MM_{X_1} \map{V} \MM_X $$
Since $i^*(x_{00}:x_{01}:x_{10}:x_{11}) = (x^2_{00}:x^2_{01}
:x^2_{10}:x^2_{11})$, we see that in terms of the canonical
Theta coordinates
$$F^*_a(x) = (\ldots,\lambda_g P^2_g(x), \ldots ).$$
Let $k_0$ be the subfield of $k$ generated by the
constants $(\lambda_g)_{g \in G}$ and $k_0^{(n)}$, for $n \geq 1$,
be the finite field extension of $k_0$ generated by the
coordinates $(x_g)_{g \in G}$ such that $F_a^{*(n)}(x) = (1:1:1:1)$
with $x = (x_g)_{g \in G}$. We obviously have a tower of extensions
$$k_0 = k_0^{(1)} \subset k_0^{(2)} \subset \cdots \subset
k_0^{(n-1)} \subset k_0^{(n)} \subset \cdots $$
and, by the computations carried out in the proof of Proposition
\ref{computations}(4), we see that $\deg [k_0^{(n)}:k_0^{(n-1)}]$
is a power of $2$. Hence, by induction, any element $x \in \pp^3 \setminus
\Omega^{Frob}$ has coodinates $(x_g)_{g \in G}$, which lie
in an extension of $k_0$ of degree $2^m$ for some $m$. But
elements of odd degree over $k_0$ are evidently dense in $\pp^3$.
\end{proof}

\begin{que}
Is $\Omega^{Rep}$ Zariski dense?
\end{que}

\section{list of questions}

\begin{enumerate}
\item
For higher genus curves we no longer have a simple description
of $\MM_X$ as for $g = 2$ (Proposition \ref{isomxp3}).
But we can ask whether the diagram
\begin{equation} \label{commdiagmod}
\begin{CD}
\MM_{X_1} @>V>> \MM_X \\
@VDVV  @VVDV \\
|2\Theta_1| @>\tilde{V}>> |2\Theta|
\end{CD}
\end{equation}
where $\tilde{V}$ is defined as in Proposition \ref{propfactorization},
is commutative for $X$ general. Note that the
``restriction'' of $V$ to $J_1$ commutes. Do we have
$\dim H^0(\MM_X, \LL^2) = 2^{g-1}(2^g +1)$? Indeed, we could check that
the latter equality implies commutativity of
\eqref{commdiagmod}.

\item
A straightforward computation shows that, for any genus $g$,
the map $\tilde{V}$ defined by the $2^g$ quadrics \eqref{polyfact}
surjects on the complement of the hyperplane $H: \ x_0 = 0$.
A priori this is not sufficient to deduce that the map
$V: \MM_{X_1} \ra \MM_X$ surjects on the complement of the
divisor $\tilde{\Theta} = D^*(H) \subset \MM_X$ as for
$g=2$. We optimistically conjecture
\begin{conj}
For any semi-stable bundle $E \in \MM_X$ satisfying
$h^0(X,E \otimes B) =0$, there exists a semi-stable bundle
$E_1 \in \MM_{X_1}$ such that $F^* E_1 = E$.
\end{conj}

As in the proof of Proposition \ref{surjectivity}, this conjecture
implies that $V : \NN_{X_1} \ra \NN_X$ is surjective.

\item
What happens for non-ordinary curves?
\end{enumerate}
We plan to return to these questions in a future work.

\bigskip

\flushleft{Yves Laszlo \\
Universit\'e Paris-Sud \\
Math\'ematiques B\^atiment 425 \\
91405 Orsay Cedex France \\
e-mail: Yves.Laszlo@math.u-psud.fr}

\bigskip

\flushleft{Christian Pauly \\
Laboratoire J.-A. Dieudonn\'e \\
Universit\'e de Nice Sophia Antipolis \\
Parc Valrose \\
06108 Nice Cedex 02 France \\
e-mail: pauly@math.unice.fr }

\end{document}